\documentclass{amsart}
\usepackage{amssymb,amsmath,amsthm,graphics,amscd,amsfonts}

\usepackage{color}

\theoremstyle{plain}
\newtheorem{theorem}{Theorem}[section]
\newtheorem{proposition}[theorem]{Proposition}

\newtheorem{lemma}[theorem]{Lemma}

\newcommand{\bfo}{{\bf o}}

\newcommand{\bfC}{{\mathbb C}}

\newcommand{\bfR}{{\mathbb R}}
\newcommand{\bfZ}{{\mathbb Z}}

\newcommand{\barg}{{\overline g}}

\newcommand{\barj}{{\overline j}}

\newcommand{\barl}{{\overline \ell}}

\newcommand{\barz}{{\overline z}}

\newcommand{\txi}{{\widetilde \xi}}

\newcommand{\mapright}[1]{\smash{\mathop{   \hbox to 0.7cm{\rightarrowfill}}
  \limits^{#1}}}

\def\grad{\mathrm{grad}}
\def\Ric{\mathrm{Ric}}

\def\a{\alpha}

\def\b{\beta}
\def\bp{\overline{\partial}}

\def\CP{{\mathbb C \mathbb P}}

\def\l{\lambda}

\def\o{\omega}

\def\om{\omega}

\def\p{\partial}
\def\bp{\overline{\partial}}

\def\s{\sigma}

\def\vph{\varphi}

\def\wt{\widetilde}

\def\calL{\mathcal L}

\begin{document}

\title{Constructing K\"ahler-Ricci solitons from Sasaki-Einstein manifolds}
\author{Akito Futaki}
\address{Department of Mathematics, Tokyo Institute of Technology, 2-12-1,
O-okayama, Meguro, Tokyo 152-8551, Japan}
\email{futaki@math.titech.ac.jp}
\author{Mu-Tao Wang}
\address{Department of Mathematics, Columbia University,
2990 Broadway
New York , NY 10027 }
\email{mtwang@math.columbia.edu }

\date{October 14, 2009}

\begin{abstract}
We construct gradient K\"ahler-Ricci solitons on Ricci-flat K\"ahler cone manifolds and
on line bundles over toric Fano manifolds.
Certain shrinking and expanding solitons are pasted together to form eternal solutions of the Ricci flow.
The method we employ is the Calabi ansatz over Sasaki-Einstein manifolds, and the results generalize constructions
of Cao and Feldman-Ilmanen-Knopf.
\end{abstract}

\keywords{Ricci soliton, Sasaki-Einstein manifold, toric Fano manifold}

\subjclass{Primary 53C55, Secondary 53C21, 55N91 }

\maketitle

\section{Introduction}
K\"ahler-Ricci solitons are self-similar solutions of the K\"ahler-Ricci flow. They are classified as expanding, steady, and shrinking solitons for obvious reasons. The convention is that an
expanding soliton lives on $(0, \infty)$ and a shrinking soliton lives on $(-\infty, 0)$.
The self-similarity reduces the Ricci flow equation to an elliptic system for a pair $(g, X)$ consists of a K\"ahler metric $g$ and a vector field
$X$ on a background manifold. In particular, any K\"ahler-Einstein metric is a steady soliton with $X=0$. K\"ahler-Ricci
solitons arise as parabolic blow-up limits of the K\"ahler-Ricci flow near a singularity.  We refer to \cite{cao06} and \cite{FIK} and for surveys of results on K\"ahler-Ricci solitons and the role they play in the singularity study of the flow.

In this article, we construct new K\"ahler-Ricci solitons from
Sasaki-Einstein manifolds. Sasaki-Einstein manifolds are links of Ricci-flat K\"ahler cones and singularity models in Calabi-Yau
manifolds. We first show that there is an expanding soliton flowing out of the K\"ahler cone over any Sasaki-Einstein manifold. The method we
employ is an ansatz of Calabi in his construction of K\"ahler-Einstein metrics (\cite{calabi79}, see also \cite{Hwang-Singer}). This Ansatz is then applied
to circle bundles over toric Fano varieties on which possibly irregular Sasaki-Einstein metric exist by \cite{FOW}. We obtain both expanding and shrinking solitons depending on the degree of the associated
line bundle. Certain pair of shrinking and expanding solitons can be pasted together to form an eternal solution of the K\"ahler-Ricci flow which lives on $(-\infty, \infty)$ with singularities
along the zero section of the line bundle, but the shrinking solitons extend smoothly to the zero section when the Sasaki-Einstein structure is regular. These results
generalize the constructions of \cite{cao96} and \cite{FIK}. While the examples in \cite{cao96} and \cite{FIK} are rotationally symmetric, our examples in general do not carry any continuous symmetry.

Similar constructions for eternal solutions of Lagrangian mean curvature flows were discovered in \cite{jlt} and \cite{lwa}. They were shown to satisfy the Brakke flow-a weak formation for the mean curvature flow. As a weak formulation of the Ricci flow has not yet been established, we may ask to what extent our examples would
qualify as generalized solutions of the Ricci flow (see also the discussion in \S 1.2 of \cite{FIK}). It will be a good indication if the flow satisfies Perelman's monotonicity formula \cite{pe} across the singularity. This will be pursued later. At this moment, we note that the Gaussian density of Perelman's functional of known Ricci solitons are computed in \cite{chi} and \cite{ha}.

The main theorems in this paper are stated as follows.
 \begin{theorem}\label{mainthm1}
 Let $M$ be a Fano manifold of dimension $m$, and $L \to M$ be a positive line bundle with $K_M = L^{-p}$, $p \in \bfZ^+$.
 For $0 < k < p$,
 let $S$ be the $U(1)$-bundle associated with $L^{-k}$, which is a regular Sasaki manifold.
 Let $Z$ be the zero section of $L^{-k}$.
 Suppose that
 $S$ admits a possibly irregular Sasaki-Einstein metric. Then there exist
 shrinking and expanding solitons on $L^{-k} - Z$, and they can be pasted together to form
  an eternal solution of the K\"ahler-Ricci flow on
 $(L^{-k} - Z) \times (-\infty, \infty)$. If $S$ admits a regular Sasaki-Einstein metric, i.e. if $M$ admits a K\"ahler-Einstein metric
 then the solution of the K\"ahler-Ricci flow corresponding to the shrinking soliton for $t \in (-\infty, 0)$ extends smoothly to the zero section $Z$.
 \end{theorem}
 Recall that a Sasaki manifold $S$ is an odd dimensional Riemannian manifold with its cone $C(S)$ a K\"ahler manifold. In the following theorem the apex of $C(S)$
 is not included in $C(S)$.
 \begin{theorem}\label{mainthm2}  Let $S$ be a compact Sasaki manfiold such that the transverse K\"ahler metric $g^T$ satisfies
Einstein equation
$$ \mathrm{Ric}^T = \kappa g^T$$
for some $\kappa < 0$ where $\mathrm{Ric}^T$ denotes the transverse Ricci curvature. Then there exists a complete
expanding soliton on the K\"ahler cone $C(S)$.
 \end{theorem}
 Note that, when a Sasaki manifold $S$ satisfies the assumption of Theorem \ref{mainthm2}, $S$ is necessarily quasi-regular so that $S$ is an
orbi-$U(1)$-bundle over an K\"ahler-Einstein orbifold with negative scalar curvature.

This paper is organized as follows. In section 2 we review K\"ahler-Ricci flows and K\"ahler-Ricci solitons.
In section 3 we review Sasaki manifolds with transverse K\"ahler-Einstein structure. In section 4 we obtain an ordinary differential equation to
get a gradient K\"ahler-Ricci soliton by Calabi's ansatz. In section 5 we extend Cao's work \cite{cao96} to construct expanding solitons
on Ricci-flat K\"ahler cones. After preparatory arguments in the case of
line bundles over Fano manifolds in section 6, we extend in section 7 the results of \cite{FIK} to construct shrinking and expanding solitons on line bundles
over Fano manifolds such that the associated $U(1)$-bundles admit Sasaki-Einstein metrics. This last condition is satisfied when the base manifolds are
toric Fano manifolds (\cite{FOW}). The shrinking soliton in section 7 and the expanding soliton in section 5 are pasted together to give an eternal solution,
and obtain Theorem \ref{mainthm1}.
In section 8 we introduce {\it gradient scalar solitons} and set up an ordinary differential equation to obtain them by Calabi's ansatz. We get a necessary
condition to have a complete gradient scalar soliton on the cone $C(S)$ of a Sasaki manifold $S$ with transverse K\"ahler-Einstein metric. We show that
a special case when the transverse K\"ahler-Einstein metric has negative transverse Ricci
 curvature gives gradient expanding Ricci solitons in Theorem \ref{mainthm2}.

\section{K\"ahler-Ricci flows and K\"ahler-Ricci solitons}
Given a K\"ahler manifold,
the K\"ahler metric $g$ can be written as
$$ g_{i\barj} = g(\frac \p{\p z^i}, \frac \p {\p \overline{z^j}})$$
where $z^1, \cdots, z^m$ are local holomorphic coordinates.
The K\"ahler form $\omega$ of $g$ is written as
$$ \omega = i g_{i\barj} dz^i \wedge d\overline{z^j},$$
and the Ricci form $\rho(\omega)$ of $\omega$ is expressed as
$$ \rho(\omega) = - i\p\bp \log \det (g_{i\barj}).$$
The coefficients
$$ R_{i\barj} = - \frac{\p^2}{\p z^i\b \overline{z^j}} \log \det (g_{k\barl}) $$
of $\rho(\omega) $ constitute the Ricci tensor $ \Ric_g$ of $g$.

A K\"ahler-Ricci flow is a family  $\om_t$ of K\"ahler forms with real parameter $t$ satisfying
\begin{equation}\label{KRF1}
\frac d{dt} \om_t = -  \frac12 \rho(\om_t).
\end{equation}
We could remove the coefficient $1/2$ in (\ref{KRF1}) by taking homothety of the K\"ahler form, but
we use the convention of (\ref{KRF1}) in order to
adapt to the convention of the paper \cite{futaki07.2}
so that we can refer to the computations there directly.

A K\"ahler-Ricci soliton is a K\"ahler form $\om$ satisfying
\begin{equation}\label{KRF2}
-\frac12 \rho(\om) = \l \om + {\mathcal L}_X \om
\end{equation}
for some holomorphic vector filed $X$ where $\l = 1, 0$ or $-1$.
Note that the imaginary part of $X$ is necessarily a Killing vector field, i.e. an infinitesimal isometry.
When
$$ \calL_X \om = i \p\bp u$$
for some real function $u$, we say that the K\"ahler-Ricci soliton is a gradient K\"ahler-Ricci soliton.
According as $\l = 1, 0$ or $-1$ the soliton is said to be expanding, steady and shrinking.

Given a K\"ahler-Ricci soliton (\ref{KRF2}) with $\l = \pm 1$, if we put
\begin{equation}\label{KRF3}
\om_t := \l t \gamma_t^{\ast} \om
\end{equation}
where $\gamma_t$ is the flow generated by the time dependent vector field
\begin{equation}\label{KRF4}
Y_t := \frac 1{\l t}X,
\end{equation}
then $\om_t $ is a K\"ahler-Ricci flow. The K\"ahler form is of course a positive form, and therefore
when $\l =1$, the Ricci flow exists for $t > 0$ and $\om_1= \om$, and
when $\l =-1$, the Ricci flow exists for $t < 0$ and $\om_{-1}= \om$.
When $\l =0$  if we put
\begin{equation}\label{KRF5}
\om_t := \gamma_t^{\ast} \om
\end{equation}
where $\gamma_t$ is the flow generated by the vector field $X$,
then $\om_t $ is a K\"ahler-Ricci flow.

\section{Ricci-flat K\"ahler cones with aperture}

In this section we first review basic facts about Sasakian geometry. Good references are the book \cite{BGbook}
and the papers \cite{FOW} and \cite{futaki07.2}.

We wish to construct gradient K\"ahler-Ricci solitons on Ricci-flat K\"ahler cones.
Recall that a cone manifold $C(S)$ is a Riemannian manifold diffeomorphic to $(0, \infty) \times S$ with a
cone metric $\barg$ is of the form
$$ \barg = dr^2 + r^2 g$$
where $g$ is a Riemannian metric on $S$ and $r$ is a coordinate on $(0, \infty)$. A Riemannian manifold $S$ is called a
Sasaki manifold if $C(S)$ is a K\"ahler cone manifold.

Let the complex dimension of $C(S)$ be $m+1$. Then
the real dimension of Sasaki manifold $(S,g)$ is $2m+1$.
$(S,g)$ is isometric to the submanifold $\{r = 1\} = \{1\} \times S
\subset (C(S),\overline{g})$, and they are usually identified.
Let $J$ be the complex structure on
$C(S)$ such that
$(C(S),\,J,\,\bar{g})$
is K\"ahler.
Then we have the vector field $\Tilde{\xi}$ and the $1$-form $\Tilde{\eta}$
on $C(S)$ defined by
$$\Tilde{\xi}=Jr\frac{\partial }{\partial r},\ \
\Tilde{\eta}=\frac{1}{r^2}\bar{g}(\Tilde{\xi},\cdot)
=\sqrt{-1}(\bar{\partial}-\partial)\log r.$$
It is easily seen that the restrictions $\xi=\Tilde{\xi}_{|S}$ and
$\eta=\Tilde{\eta}_{|S}$ to $\{r=1\}\simeq S$ give a vector field and a
 one form on $S$.
 The one form $\eta$ on $S$ is a contact form and
 the vector field $\xi$ is the Reeb vector field of the contact form $\eta$, that is $\xi$ is the unique vector field
 which satisfies
 $$ i(\xi)\eta = 1\ \qquad \mathrm{and}\qquad  i(\xi)d\eta = 0.$$

There are two K\"ahler structures involved in the study of Sasaki manifolds. One is the K\"ahler structure on $C(S)$.
The K\"ahler form $\omega$ of
$(C(S),\,J,\,\bar{g})$
is given by
$$
\omega=\frac12 d(r^2\Tilde{\eta})=\frac{\sqrt{-1}}{2}\partial
\bar{\partial}r^2.$$
The second one is the transverse K\"ahler structure of the flow, called the Reeb flow, generated by the Reeb vector field $\xi$.
This is a collection of K\"ahler structures on the local orbit spaces of the Reeb flow.
The vector field $\Tilde{\xi}$ is a Killing vector field on
$(C(S),\bar{g})$ with the length $\bar{g}(\Tilde{\xi}, \Tilde{\xi})^{1/2}=r$.
The complexification $\Tilde{\xi}-\sqrt{-1}J\Tilde{\xi}$
of the vector field is holomorphic on $(C(S),J)$. Since the local orbit spaces of the Reeb flow on $S$ and
the local orbit spaces of the holomorphic flow generated by $\Tilde{\xi}-\sqrt{-1}J\Tilde{\xi}$ on $C(S)$
along $S$ can be identified then
they define a transverse holomorphic structure of the Reeb flow, i.e.
holomorphic structures on the local orbit spaces of the Reeb flow. But since $\eta$ is a contact form and non-degenerate on
the contact distribution, i.e. the kernel of $\eta$,
we obtain the transverse K\"ahler structure by identifying the tangent spaces of local orbit spaces with the contact distribution.
Thus the transverse K\"ahler structure is a collection of K\"ahler structures on local orbit
spaces of the Reeb flow. The K\"ahler forms on local orbit spaces are lifted to $S$ to form a global $2$-form
$$ \omega^T := \frac 12 d\eta,$$
called the transverse K\"ahler form. The transverse K\"ahler form can be lifted also to $C(S)$ and can be expressed as
$$ \wt \omega^T = \frac12 d\wt{\eta} = \frac i2 d(\bp - \p) \log r = dd^c \log r.$$
The Ricci curvature $\Ric^T$ of the transverse K\"ahler metric $g^T$ is related to the Ricci curvature $\Ric_g$ of $(S,g)$ by
\begin{equation}\label{sasaki3}
\mathrm{Ric}_g = \mathrm{Ric}^{T} - 2g^{T} + 2m\ \eta\otimes \eta.
\end{equation}
We wish to use Calabi's ansatz when the transverse K\"ahler structure is K\"ahler-Einstein. This last condition
is equivalent to say that the orbit spaces have
K\"ahler-Einstein metrics.

There is a related notion in classical Sasakian geometry, called $\eta$-Einstein manifolds.
A Sasaki manifold $S$ is called an $\eta$-Einstein manifold if for some constants $\a$ and $\b$
$$ \Ric_g = \alpha g + \beta\,\eta \otimes \eta.$$
Since $\Ric(\xi,\xi) = 2m$ by (\ref{sasaki3}), we have $\a + \b = 2m$, and (\ref{sasaki3}) also shows that
the transverse K\"ahler metric is Einstein with
$$ \rho(\omega^T) = (\a + 2) \omega^T.$$
Conversely a Sasaki manifold with a transverse K\"ahler-Einstein metric has an $\eta$-Einstein metric.
Therefore $\Ric^T$ is positive definite if and only if $\alpha + 2 > 0$, that is $\alpha > -2$.
Obviously $(S,g)$ is Sasaki-Einstein, i.e.
$\eta$-Einstein with $\beta = 0$ if and only if
$\Ric = 2m g$, and also if and only if $\Ric^T = (2m+2) g^T$.
The Gauss equation also tells us that $(S, g)$ is a Sasaki-Einstein manifold if and only if $(C(S), \barg)$ is a Ricci-flat K\"ahler manifold.

The typical example is when $(S, g)$ is the $(2m+1)$-dimensional standard sphere in which case $C(S) = \bfC^{m+1} - \{\bfo\}$ with the flat metric
and the orbit space of the Reeb flow is $\CP^m$ with a multiple of the Fubini-Study metric such that the Einstein constant is $(2m+2)$.

Given a Sasaki manifold with the K\"ahler cone metric $\barg = dr^2 + r^2 g$, we transform
the Sasakian structure by deforming $r$ into $r' = r^a$ for positive constant $a$.
This transformation is called the $D$-homothetic transformation. Then the new
Sasaki structure has
\begin{equation}\label{s12}
\eta' = d \log r^a = a\eta, \quad \xi' = \frac 1a \xi,
\end{equation}
\begin{equation}\label{s13}
g' = a g^T + a\eta \otimes a\eta = ag + (a^2 -a)\eta\otimes\eta.
\end{equation}
Suppose that $g$ is $\eta$-Einstein with $\mathrm{Ric}_g = \a g + \b \eta \otimes \eta$.
Since the Ricci curvature of a K\"ahler manifold is invariant under homotheties
the transverse Ricci form is invariant under the D-homothetic transformations :
$\mathrm{Ric}^{\prime T} = \mathrm{Ric}^T$.
From this and $\mathrm{Ric}_{g'}(\xi',\xi') = 2m$ we have
\begin{eqnarray}\label{s14}
\mathrm{Ric}_{g'} &=& \mathrm{Ric}^{\prime T} - 2g^{\prime T} + 2m \eta' \otimes \eta' \\
&=& \mathrm{Ric}^T - 2a g^T + 2m \eta' \otimes \eta' \nonumber\\
&=& \mathrm{Ric}|_{D \times D} + 2g^T - 2a g^T + 2m \eta' \otimes \eta' \nonumber\\
&=& \a g^T  + 2g^T - 2a g^T + 2m \eta' \otimes \eta'. \nonumber
\end{eqnarray}
This shows that $g'$ is $\eta$-Einstein with
\begin{equation}\label{s15}
\a' =  \frac{\a + 2 - 2a}a.
\end{equation}
Thus we have proved the following well-known fact.
\begin{lemma}\label{sasaki1}
Under the $D$-homothetic transformation of an $\eta$-Einstein metric we have
a new $\eta$-Einstein metric with
\begin{equation}\label{s16}
\rho^{\prime T} = \rho^T, \quad \omega^{\prime T} = a\omega^T, \quad
\rho^{\prime T} = (\a' + 2) \omega^{\prime T} = \frac{\a + 2}a \omega^{\prime T},
\end{equation}
and thus, for any positive constants $\kappa$ and $\kappa'$, a transverse K\"ahler-Einstein metric
with Einstein constant $\kappa$ can be transformed by a $D$-homothetic transformation to a
transverse K\"ahler-Einstein metric with Einstein constant $\kappa'$. The same is true for negative $\kappa$ and $\kappa'$,
and for $\kappa = \kappa' = 0$.
\end{lemma}

 Given a Sasaki-Einstein metric, a positive constant multiple of the K\"ahler cone metric of its D-homothetic transformation may be called a {\it Ricci-flat
 K\"ahler cone metric with aperture}, whose K\"ahler form $\wt \o$ is of the form
 \begin{equation}\label{sasaki2}
  \wt\o =  C i\p\bp \log r^{2q}
  \end{equation}
 for some positive constants $C$ and  $q$ where $\barg = dr^2 + r^2 g $
 is the Ricci-flat K\"ahler cone metric of the given Sasaki-Einstein metric.

\section{Solitons on Ricci-flat K\"ahler cones}\label{SRFKC}

In this section we apply Calabi's ansatz to transversely K\"ahler-Einstein Sasaki manifolds to obtain gradient K\"ahler Ricci solitons.
Suppose that we have a  transversely K\"ahler-Einstein Sasaki manifold so that we have an $\eta$-Einstein metric with
$$ \Ric_g = \alpha g + \beta\,\eta \otimes \eta.$$
This is a transverse K\"ahler-Einstein
metric with transverse K\"ahler form
$\omega^T = \frac12 d\eta$ satisfying
$$\rho^T = \kappa\ \omega^T$$
with $\kappa = \alpha + 2$.
As was explained in the previous section,
on $C(S)$ we have
$$ \eta = 2 d^c \log r, \qquad \omega^T = dd^c \log r. $$
The Calabi's ansatz seeks for a special metric of the form
$$\omega = \om^T + i\p\bp F(s).$$
where we put
$$s = \log r$$
and where
$$F \in C^{\infty}((s_1,s_2)), \qquad (s_1, s_2) \subset (-\infty, \infty).$$
Here we search metrics of this form for K\"ahler-Ricci solitons. We further put
$$\sigma = 1 + F'(s),$$
and define $\varphi(\s)$ by
\begin{equation}\label{sasaki4}
\varphi(\sigma) = F''(s).
\end{equation}
Since
$$  i\p\bp F(s) = i F''(s) \p s \wedge \bp s + i F'(s) \p\bp s $$
we have
$$ \omega = \sigma \o^T + \varphi(\sigma) \p s \wedge \bp s.$$
Put
$$\lim_{s \to s_1} \s(s) = \lim_{s \to s_1} (1+F'(s)) = a, \qquad \lim_{s \to s_2} \s(s) = \lim_{s \to s_2} (1+F'(s)) = b.$$
Because of the positivity of $\omega$ we must have
\begin{equation}\label{Cone0}
 \s > 0 \qquad \mathrm{and} \qquad \vph(\s) > 0
 \end{equation}
on the region $a < \s < b$.
By $\sigma'(s) > 0$, the map $\s : (s_1, s_2) \to (a,b)$ is a diffeomorphism.

Conversely, given a positive function $\varphi$ on $(a,b)$ with $a > 0$ such that
$$ \lim_{\sigma \to a^+} \int_{\sigma_0}^\sigma \frac{dx}{\varphi(\s)} = s_1, \quad \lim_{\sigma \to b^-} \int_{\sigma_0}^\sigma \frac{dx}{\varphi(\s)} = s_2$$
we define $\sigma(s)$ by
$$ s = \int_{\s_0}^{\s(s)} \frac {dx}{\varphi(x)}  $$
and define $F(s)$ by
$$ F(s) = \int_{\s_0}^{\s(s)} \frac {x-1}{\varphi(x)} dx.$$
Note that
$$ \frac {ds}{d\s} = \frac 1{\vph(\s)}, \qquad \frac{d\s}{ds} = \vph(\s) = F''(s) $$
$$ \frac{dF}{ds} = \frac{\s -1}{\vph(\s)} \frac{d\s}{ds} = \s -1.$$
If we put $s = \log r$ and consider $\sigma$ and $F$ as a smooth function on
$$ C(S)_{(s_1,s_2)} := \{ e^{s_1} < r < e^{s_2} \} \subset C(S)$$
then
\begin{eqnarray}
\om_{\vph} &:=& \om^T + dd^c F(s)  \label{omega1}\\
&=& \s\om^T + \vph(\s) i \p s \wedge \bp s \nonumber \\
&=& \s \om^T + \vph(\s) \om_{\mathrm{cyl}} \nonumber
\end{eqnarray}
gives a K\"ahler form on $C(S)_{(s_1,s_2)}$ with the K\"ahler metric
\begin{equation*}
g = \s g^T + \vph(\s) g_{\mathrm{cyl}}
\end{equation*}
where $\om_{cyl}$ and $g_{cyl}$ denote the cylindrical K\"ahler form and K\"ahler metric on $\bfC - \{\bfo\}$.
It is also possible to write
\begin{equation}\label{omega}
\om_{\vph} =  i\p\bp (s + F(s)).
\end{equation}
Let
$(z^0, z^1, \cdots , z^m) $ be local holomorphic coordinates on $C(S)$
such that
$$\frac \p{\p z^0} = \frac 12 (r\frac \p{\p r} - iJ r\frac \p{\p r}) = \frac 12 (r\frac \p{\p r} - i \wt\xi) .$$
Then we have
$$ dz^0 = \frac {dr}r + i \eta $$
and
$$ i dz^0 \wedge d\overline{z}^0 = 2 \frac{dr}r \wedge \eta. $$
Using these coordinates one can show as in \cite{futaki07.2} that
\begin{eqnarray}
\omega_\vph ^{m+1} &=& \s^m(m+1)\vph(\s) i\p s \wedge \bp s \wedge (\om^T)^m \nonumber \nonumber\\
&=& \s^m (m+1) \vph(\s) \frac i2 dz^0 \wedge d\barz^0 \wedge (\om^T)^m
\end{eqnarray}
and that
\begin{eqnarray}
\rho_{\vph} &=& \rho^T - i \p\bp \log (\s^m \vph(\s))\nonumber\\
&=& \kappa \om^T - i \p\bp \log (\s^m \vph(\s)) \label{ricci}.
\end{eqnarray}

\begin{lemma}\label{grad1}  Let $Q$ be a smooth function in $s$. Then
$\grad'\ Q(s)$ is a holomorphic vector field if and only if
$$ \grad\ Q(s) = \mu r \frac \p{\p r}$$
for some constant $\mu \in \bfR$. Moreover this is equivalent to
\begin{equation}\label{grad2}
Q_s = \mu \vph(\s)
\end{equation}
and also to
\begin{equation}\label{grad3}
Q = \mu \s + c
\end{equation}
where $c$ is a constant.
\end{lemma}

\begin{proof} The former half follows because $r\frac \p{\p r} - iJ r\frac \p{\p r}$ is the only $s$-dependent holomorphic vector field.
For the latter half of the lemma,
since
\begin{eqnarray*}
\vph(\s) g_{\mathrm{cyl}} &=& \vph(\s) (ds^2 + d\theta^2 ) \\
&=&  \vph(\s) \left( \frac {dr^2}{r^2} + d\theta^2\right) = \frac{\vph(\s)}{r^2} (dr^2 + r^2 d\theta^2)
\end{eqnarray*}
we have
$$ \grad\,Q = \frac {r^2}{\vph(\s)} Q_s \frac{\p s}{\p r} \frac \p{\p r} = \frac{r^2}{\vph(\s)} Q_s \frac 1r \frac \p{\p r}  = \frac r{\vph(\s)} Q_s \frac \p{\p r}. $$
This and the former half of the lemma show
$$ \frac r{\vph(\s)} Q_s = \mu r,$$
which implies (\ref{grad2}). This completes the proof of Lemma \ref{grad1}.
\end{proof}
Now let us derive the gradient K\"ahler-Ricci soliton equation in terms of $\vph$.
By (\ref{ricci}) we have
 $$\rho_{\om_{\vph}} = i \p\bp (\kappa s - \log (\s^m \vph(\s))).$$
 Comparing this with (\ref{omega})
\begin{equation}\label{ode2}
\rho_{\omega_{\vph}} + 2\lambda \omega_{\vph} = i\p\bp (\kappa s - \log (\s^m \vph(\s)) + 2\lambda s + 2\lambda F).
\end{equation}
Put
$$ Q := - \kappa s + \log (\s^m \vph(\s)) - 2\lambda s - 2\lambda F$$
so that
$$ \rho_{\omega_{\vph}} + 2\lambda \omega_{\vph} = -i\p\bp Q.$$

In order for $\om_{\vph}$ to be a gradient K\"ahler-Ricci soliton, $\grad\, Q$ must be the real part of a holomorphic vector field.
So we may apply Lemma \ref{grad1} to this $Q$. Then we see from (\ref{grad2}) that $\vph(\sigma)$ must satisfy
\begin{equation}\label{GKRS2}
\vph'(\s) + \left( \frac m \s - \mu \right) \vph(\s) - (\kappa + 2\l \s ) = 0.
\end{equation}
In general a solution to the ODE $y' + p(x)y = q(x)$ is
\begin{equation}\label{ode}
 y = e^{- \int p(x) dx} \left( \int q(x) e^{\int p(x) dx} dx +C\right) .
 \end{equation}
It follows from (\ref{ode}) that the solution to (\ref{GKRS2}) is given by
\begin{equation}\label{GKRS3}
\vph(\s) = \frac{\nu e^{\mu\s}}{\s^m} - \frac{2\l\s}{\mu} - \frac{2\l + \frac{\kappa\mu}{m+1}}{\mu^{m+2}} \sum_{j=0}^m \frac{(m+1)!}{j!} \mu^j \s^{j-m}.
\end{equation}

In the next section we construct K\"ahler-Ricci solitons on $C(S)$. We thus assume $(s_1, s_2) = (-\infty. \infty)$.
We also assume $a := \lim_{s \to -\infty}\s(s) \ge 0$. It follows from this assumption that
\begin{equation}\label{Cone1}
\vph(a) = 0.
\end{equation}
It also follows from (\ref{GKRS3}) and (\ref{Cone1}) that, in either cases of $\l = -1, 0$ or $\l = 1$, $\nu$ is determined by
\begin{eqnarray}\label{GKRS5}
 \nu &=& e^{-\mu a} a^m \left( \frac {2\l}\mu a + \frac{2\l +\frac{\kappa\mu}{m+1}}{\mu^{m+2}} \sum_{j=0}^m \frac{(m+1)!}{j!} \mu^j a^{j-m}\right)\\
 &=:& \nu_a^\l(\mu).\nonumber
 \end{eqnarray}

\section{Expanding solitons on Ricci-flat K\"ahler cones}

 In this section we extend Cao's construction \cite{cao96} of expanding soliton on $\bfC^n$ to the general Ricci-flat K\"ahler cones, i.e. the K\"ahler cones over
 Sasaki-Einstein manifolds.
 In this case we require
 $$a = 0.$$
 Geometric reasoning of this requirement is given in sections 4.1 and 4.2 in \cite{FIK}, and we do not reproduce it here.
 Then near $\s = 0$ we have
 \begin{eqnarray}
\vph(\s) &=& \frac{\nu e^{\mu\s}}{\s^m} - \frac{\nu_0^\l(\mu)}{\s^m} \sum_{j=0}^m \frac{\mu^j \s^j }{j!} - \frac{2\l\s}{\mu} \label{Cone2}\\
&\approx& \frac{\nu - \nu_0^\l(\mu)}{\s^m} \label{Cone3}
\end{eqnarray}
where $\nu_0^\l(\mu)$ is the one given as (\ref{GKRS5}) with $a = 0$:
$$ \nu_0^\l(\mu) = \frac{(m+1)!(2\l + \frac{\kappa\mu}{(m+1)})}{\mu^{m+2}}.$$
But we must have
\begin{equation}\label{Cone4}
\nu = \nu_0^\l(\mu),
\end{equation}
for if $\nu < \nu_0^\l(\mu)$ then $\vph(\s) < 0$ near $\s = 0$ contradicting (\ref{Cone1}),
and if $\nu > \nu_0^\l(\mu)$ then $\s(s)$ becomes $0$ for finite $s > -\infty$.

With this $\nu$, $\vph$ is written as
$$ \vph(\s) = \frac \kappa{m+1} \s + \nu_0^\l(\mu) \sum_{j=2}^\infty \frac{\mu^{j+m}\s^j}{(j+m)!},$$
so we have a solution $\s : (-\infty, r_1] \to \bfR$ of $d\s/ds = \vph(\s)$.

Put $\l = 1$ and $\mu = -  1/q$ with
arbitrary fixed $q > 0$.
We know that $\vph(0) = 0$ and $\vph'(0) > 0$.
If $b$ is any positive solution of $\vph(b) = 0$ then by (\ref{GKRS2})
$$ \vph'(b) = \kappa + 2b > 0.$$
This implies that there is no positive zero $b$, and we have $\vph > 0$ for all $\s > 0$. For large $\s$ we see from (\ref{Cone2}) that
$$ \frac{d\s}{ds} = - \frac{2\s}{\mu} + G(\frac 1\s)$$
where $G$ is smooth at zero.
Since $\mu < 0$ then $\s$ extends for all $s \in \bfR$ and has the form
$$ \s(s) = e^{-\frac {2s}\mu}E(e^{\frac {2s}\mu}) = r^{-\frac 2\mu} E(r^{\frac2\mu})$$
where $E$ is smooth at zero and $E(0) > 0$. This soliton is asymptotic to a Ricci-flat K\"ahler cone metric with aperture as can be seen as follows.
The vector field $\frac{\mu}{2t} r\frac \p{\p r}$ generates a one parameter group $\{\gamma_t\}$ of transformations such that
$$ \gamma_t^{\ast} r = t^\frac\mu2 r.$$
Then
$$ \gamma_t^{\ast} s = \log \gamma_t^{\ast} r = \frac\mu2 \log t + s.$$
The flow $\{t\gamma_t^{\ast}\omega\}$ satisfies
\begin{eqnarray}
t\gamma_t^{\ast} \omega &=& t (t^{\frac\mu2}r)^{-\frac 2\mu} E((t^\frac\mu2 r)^{\frac 2\mu}) i\p\bp s  \nonumber\\
&& + t \left(-\frac 2\mu t^{-1}r^{-\frac 2\mu} E(tr^{\frac2\mu}) - \frac 4{\mu^2} E'(tr^{\frac2\mu})\right) i\p s \wedge \bp s \nonumber \\
&\to & E(0)\left( r^{-\frac 2\mu} i \p\bp s - \frac 2\mu r^{- \frac 2\mu} i\p s \wedge \bp s \right) \label{limit1}
\end{eqnarray}
as $t \to 0$.
 Since we put $q = - \frac 1{\mu}$ then (\ref{limit1}) is equal to
 \begin{equation}
 E(0) \left(r^{2q} i \p\bp \log r + 2q r^{2q} i\p \log r \wedge i\bp \log r \right)
 = E(0) i\p\bp \left(\frac {r^{2q}}{2q}\right).
\end{equation}
This is a K\"ahler cone metric of an $\eta$-Einstein Sasaki manifold.

Thus we have proved the following:
\begin{theorem}\label{expandingcao}
Let $S$ be a compact Sasaki-Einstein manifold and $C(S)$ its K\"ahler cone. Then there is a gradient expanding soliton
which is asymptotic to a Ricci-flat cone metric with aperture.
\end{theorem}

 \section{Line bundles over toric Fano manifolds}\label{linebundle}

 Let $M$ be a Fano manifold of dimension $m$, and $L \to M$ be a positive line bundle with $K_M = L^{-p}$, $p \in \bfZ^+$.
 Take $k \in \bfZ^+$. Let $S$ be the $U(1)$-bundle associated with $L^{-k}$, which is a regular Sasaki manifold with the K\"ahler cone $C(S)$ biholomorphic to $L^{-k}$ minus
 the zero section. It is proven in \cite{FOW} that when $M$ is toric then  $S$ admits a possibly irregular toric Sasaki-Einstein metric.
 Keeping this result in mind we assume that
 $S$ is a possibly irregular Sasaki-Einstein manifold whose cone $C(S)$ is nevertheless biholomorphic
 to the cone of the regular Sasaki structure, i.e.
 the total space of $L^{-k}$ minus zero section.

 Let $\kappa = \frac{2p}k$. Then by a $D$-homothetic transformation we may assume we have a transverse K\"ahler-Einstein metric, i.e. $\eta$-Einstein
 Sasaki metric, such that
 $$ \rho^T = \kappa \omega^T$$
 where $\omega^T$ and $\rho^T$ are respectively the transverse K\"ahler form and its transverse Ricci form. Then we have
 \begin{equation}\label{omega^T}
  2[\om^T] = c_1(L^k).
  \end{equation}
 In this set-up we apply the computations in section 4, and we have a gradient K\"ahler-Ricci soliton $\om_{\vph}$ with $\vph$ given by (\ref{GKRS3}).
 Let
 $$ \lim_{s\to -\infty} \s(s) = a, \qquad \lim_{s \to \infty} \sigma(s) = b, $$
 and suppose that $a>0$.
 \begin{lemma}\label{zero}
If $\s$ is a zero of $\varphi$ then
\begin{equation}\label{GKRS4}
\vph'(\s) = \kappa + 2\l\s.
\end{equation}
Thus if $\l = -1$ then there are at most two positive zeros, one $a$ with $0 < a \le \frac\kappa2$ and one $b$ with $\frac\kappa2 \le b$.
If $\l = 1$ then there is at most one positive zero $0 < a < 1$.
\end{lemma}
\begin{proof} This follows immediately from (\ref{GKRS2}). \end{proof}
 \begin{theorem}\label{extension}  Suppose that the Sasaki-Einstein structure is regular. Then
the K\"ahler-Ricci soliton given by the solution of (\ref{GKRS3}) with $a>0$ extends to the zero section smoothly if and only if
$$ a = \l(1 - \frac pk ) = \l(1 - \frac \kappa 2).$$
In particular we have  $0 < k < p$ if $\l = -1$ and that $p < k $ if $\l =1$.
\end{theorem}
 \begin{proof}
Suppose that $\om_{\vph}$ extends to the zero section as a K\"ahler form. Since $\om_{\vph}$ satisfies the K\"ahler-Ricci soliton equation
we have
$$ - [2\l \om_{\vph}]|_M = [\rho_{\om_\vph}]|_M = c_1(M) + c_1(L^{-k}).$$
On the other hand, from
$$ \omega_{\vph} = \s \omega^T + \varphi(\s) \omega_{\mathrm{cyl}} $$
and  (\ref{omega^T}) we have
$$ - [2\l \om_{\vph}]|_M = - \l ak c_1(L).$$
If $\l = \pm 1$ then  we have
$$ a = \l(1 - \frac pk ) = \l(1 - \frac \kappa 2).$$
Since $a > 0$ this shows that $0 < k < p$ if $\l = -1$ and that $p < k $ if $\l =1$.

Conversely, suppose that we have $ a = \l(1 - \frac \kappa 2)$.
By (\ref{GKRS4}) if $\l = -1$ we have
$$\vph'(a) = \kappa - 2a = \kappa - (\kappa -2) = 2. $$
If $ \l = 1$ then
$$\vph'(a) = \kappa + 2a = \kappa + (2 - \kappa) = 2. $$
Then the extension to the zero section follows from the Proposition \ref{C30} below.
\end{proof}

Here we take up the problem of completeness of the metrics obtained by Calabi's ansatz
starting from a compact $\eta$-Einstein metric. We do not necessarily assume that the
$\eta$-Einstein structure has transversely positive Ricci curvature.
\begin{proposition}\label{C29}Let $\omega_{\varphi}$ be the K\"ahler form obtained by
Calabi's ansatz starting from a compact Sasaki manifold with
an $\eta$-Einstein metric $g$. Then $\omega_{\varphi}$ defines a
complete metric with noncompact ends towards the end points of $I = (a,b)$
if and only if the following conditions are satisfied at the end points:
\begin{itemize}
\item
At $\s = a$, $\varphi$ vanishes at least to the second order.
\item
If $b$ is finite then as $\varphi$ vanishes at $\s = b$ at least to the second order.
\item
If $b = \infty$ then $\varphi$ grows at most quadratically as $\s \to \infty$.
\end{itemize}
\end{proposition}
\noindent
Before starting the proof we change the variable by
$$ \tau = \s - a $$
because this makes the arguments more transparent. We regard $\vph$ as a function of $\tau$, and then what
we need to show is that, for example, at $\tau = 0$, $\vph$ vanishes at least to the second order, and if $b = \infty$
then $\varphi$ grows at most quadratically as $\tau \to \infty$.
\begin{proof} First define the function $\ell(s)$ by
\begin{equation}\label{C27}
\ell(s) = \int_{\tau_0}^{\tau(s)} \frac{dx}{\sqrt{\varphi(x)}}.
\end{equation}
Then
\begin{equation}\label{C28}
\frac{d\ell}{ds} = \frac{1}{\sqrt{\varphi(\tau)}}\frac{d\tau}{ds} = \sqrt{\varphi(\tau)}.
\end{equation}
Thus $\ell(s)$ gives the geodesic length along the $s$-direction with respect to the K\"ahler form $\omega_{\varphi}$ of
(\ref{omega1}); recall $s = \log r$.

Next consider at $\tau = 0$. By elementary calculus $\ell(s) \to \infty$ as $\tau \to 0$
if and only if $\varphi$ vanishes at $0$ at least to the second order. By the same reason, if
$b$ is finite then $\varphi$ must vanish at $\tau = b-a$ at least to the second order.
Similarly if $b = \infty$,  $\ell(s) \to \infty$ as $\tau \to \infty$ if and only if $\varphi$ grows at most quadratically.
\end{proof}

\begin{proposition}\label{C30}Let $\omega_{\varphi}$ be the K\"ahler form obtained by
Calabi's ansatz starting from a regular compact $\eta$-Einstein Sasaki manifold.
Suppose that the profile $\varphi$ is defined on $(a,\infty)$ and that $t_1 = - \infty$.
Then $\omega_{\varphi}$ defines a
complete metric, has a noncompact end towards $\s = \infty$
and extends to a smooth metric on the total space of the line bundle up to the zero section
if and only if $\varphi$ grows at most quadratically as $\s \to \infty$ and
$\varphi(a) = 0$ and $\varphi'(a) = 2$.
\end{proposition}
\begin{proof} As before we use the change of variable
$$ \tau = \s -a$$
and regard $\vph$ as a function of $\tau$.
As in the proof of the previous proposition $\varphi$ must grow at most
quadratically as $\tau \to \infty$. Now let us consider at $\tau = 0$.
From the assumptions of the proposition the Sasaki
manifold $S$ is the total space of the $U(1)$-bundle associated with an Hermitian line bundle $(L,h)$
whose curvature form $\omega^T$ is K\"ahler-Einstein on the base manifold of $L$.
Let $z$ be the fiber coordinate and put $r^2 = h|z|^2$. This is the K\"ahler form of the cone $C(S)$, which is
isomorphic to $L$ minus
the zero section.
Recall that
the K\"ahler form
$\omega_{\varphi}$ given by (\ref{omega1})
is of the form
$$ \omega_{\vph} = (\tau + a)\om^T + \vph(\tau)\p s \wedge \bp s.$$
Let $\pi : L^k \to M$ be the projection and $i : M \to L^k$ the inclusion to the zero section.
Since $\om^T$ is the restriction to $L^k - i(M)$ of $\pi^{\ast}\omega_{KE}$ where $\omega_{KE}$ is the K\"ahler form of a K\"ahler-Einstein metric
on the base manifold $M$, $\omega^T$ naturally extends to the zero section $i(M)$.
Therefore we have only to consider $\om_{\vph}$ in the direction
of the holomorphic flow generated by $\frac12 (\txi - iJ\txi)$.  Thus we look at
\begin{eqnarray}\label{C31}
\varphi(\tau)\, idt\wedge d^ct &=& \frac{\varphi(\tau)}4 i d\log r^2 \wedge d^c \log r^2\\
&=& \frac{\varphi(\tau)}{4r^2}  (h\,idz \wedge d\barz + O(|z|)).\nonumber
\end{eqnarray}
Hence we need only to find the condition for
$\lim_{\tau \to 0} \varphi(\tau)/r^2 $
to exist and be positive. Suppose that
\begin{equation}\label{C32}
\varphi(\tau) = a_1\tau + O(\tau^2).
\end{equation}
Since $s = \log r$, $\tau +a = 1 + F'(s)$ and $\varphi(\tau) = F''(s)$ we have
\begin{equation}\label{C33}
\frac{d\tau}{ds} =\varphi(\tau) = a_1\tau + O(\tau^2).
\end{equation}
Thus
\begin{equation}\label{C134}
\lim_{\tau \to 0}\frac{\varphi(\tau)}{r^2} =
\lim_{s \to -\infty}\frac{\varphi'(\tau)\frac{d\tau}{ds}}{2r\frac{dr}{ds}}
= \frac{a_1}2 \lim_{\tau \to 0}\frac{\varphi(\tau)}{r^2}.
\end{equation}
Therefore if $\lim_{\tau \to 0}\frac{\varphi(\tau)}{r^2}$ exists and is positive then $a_1 = 2$,
i.e. $\varphi'(0) = 2$.
Conversely if $\varphi'(0) = 2$ then we have
\begin{equation}\label{C135}
\frac{d\tau}{dt} = \varphi(\tau) = 2\tau + O(\tau^2) = 2\tau\alpha(\tau)
\end{equation}
where $\alpha(\tau)$ is a function of $\tau$ real analytic near $\tau = 0$ with $\alpha(0) = 1$
since $\varphi$ is a real analytic function by (\ref{GKRS3}).
We then have
\begin{equation}\label{C136}
\frac{d\tau}{\tau\alpha(\tau)} = 2dt
\end{equation}
and from this
\begin{equation}\label{C137}
\log \tau + \beta(\tau) = c + 2t
\end{equation}
where $\beta(\tau)$ is a real analytic function of $\tau$ with $\beta(0) = 0$. From this we have
\begin{equation}\label{C138}
\tau = e^{-\beta(\tau)} e^{c + 2t} = r^2e^{c - \beta(\tau)}.
\end{equation}
Thus we obtain
\begin{equation}\label{C139}
\lim_{\tau \to 0} \frac{\varphi(\tau)}{r^2} = \lim_{\tau \to 0}\frac{2\tau + O(\tau^2)}{r^2} = e^{c}.
\end{equation}
This completes the proof of Proposition \ref{C30}.
\end{proof}

 \section{Expanding and shrinking solitons on line bundles over Fano manifolds}

 In the set-up of the section \ref{linebundle}
 we first consider expanding solitons on $L^{-k}$ with $k > p$. Recall that $\lambda = 1$ and $a = 1 - p/k$ in this case by Theorem \ref{extension}.
 By the same argument as in the proof of Lemma 5.1 of \cite{FIK} we can prove that
an expanding soliton on $L^{-k}$ must have $\mu < 0$.
By (\ref{GKRS3}) the dominant term of $\vph$ is $- \frac {2\s} \mu$ because the exponential term is tame. We may write
$\s_s = \vph(\s)$ in the form
$$ \s_s = - \frac {2\s} \mu + G\left(\frac 1\s\right) $$
where $G$ is smooth at zero.
Considering the behavior of $1/\s$ for large $s$ we find that
$$ \s(s) = e^{-\frac {2s}\mu} B(e^{\frac {2s}\mu}) $$
where $B$ is a smooth function with $B(0) > 0$.
Using $s = \log r$ we get
\begin{eqnarray}
\omega &=& e^{-\frac {2s}\mu} B(e^{\frac {2s}\mu}) \omega^T + \left( - \frac2\mu e^{-\frac {2s}\mu}B(e^{\frac {2s}\mu}) - \frac 4{\mu^2} B'(e^{\frac {2s}\mu})\right) \omega_{\mathrm{cyl}}\\
&=& r^{-\frac 2\mu} B(r^{\frac 2\mu}) i\p\bp s + \left(-\frac 2\mu r^{-\frac 2\mu} B(r^{\frac 2\mu}) - \frac 4{\mu^2} B'(r^{\frac 2\mu})\right) i\p s\wedge \bp s.
\end{eqnarray}

The vector field $\frac{\mu}{2t} r\frac \p{\p r}$ generates a one parameter group $\{\gamma_t\}$ of transformations such that
$$ \gamma_t^{\ast} r = t^{\frac\mu2} r.$$
Then
$$ \gamma_t^{\ast} s = \log \gamma_t^{\ast} r = \frac\mu2 \log t + s.$$
The flow $\{t\gamma_t^{\ast}\omega\}$ satisfies
\begin{eqnarray}
t\gamma_t^{\ast} \omega &=& t (t^{\frac\mu2}r)^{-\frac 2\mu} B((t^\frac\mu2 r)^{\frac 2\mu}) i\p\bp s  \nonumber\\
&& + t \left(-\frac 2\mu t^{-1}r^{-\frac 2\mu} B(tr^{\frac2\mu}) - \frac 4{\mu^2} B'(tr^{\frac2\mu})\right) i\p s \wedge \bp s \nonumber \\
&\to & B(0)\left( r^{-\frac 2\mu} i \p\bp s - \frac 2\mu r^{- \frac 2\mu} i\p s \wedge \bp s \right) \label{limit2}
\end{eqnarray}
as $t \to 0$.
 If we put $q = - \frac 1{\mu}$ then (\ref{limit2}) is equal to
 \begin{equation}
 B(0) \left(r^{2q} i \p\bp \log r + 2q r^{2q} i\p \log r \wedge i\bp \log r \right)
 = B(0) i\p\bp \left(\frac {r^{2q}}{2q}\right).
\end{equation}
This is a Ricci-flat K\"ahler cone metric with aperture.

Thus we have proved\footnote{ Professor X.H.~Zhu kindly informed us that the regular case of Theorem \ref{thmexpand} and \ref{shrinking(t<0)} was also
obtained by Y. Bo \cite{bo08}}:
\begin{theorem}\label{thmexpand} Let $M$ be a Fano manifold, and $L \to M$ be a positive line bundle with
$L^{-p} = K_M$, $p \in \bfZ^+$. Suppose that the $U(1)$-bundle of $K_M$ admits a possibly irregular Sasaki-Einstein metric
whose cone $C(S)$ is biholomorphic to the total space of $K_M$ minus the zero section.
For $k > p$, $L^{-k}$ minus the zero section admits a gradient expanding soliton
such that the corresponding K\"ahler-Ricci flow $g(t)$ converges to a Ricci-flat K\"ahler cone metric with aperture,
or equivalently a K\"ahler cone metric over a transversely K\"ahler-Einstein Sasaki manifold.
Here the K\"ahler cone manifold is biholomorphic to $L^{-k}$ minus the zero section and the
transversely K\"ahler-Einstein Sasaki manifold is diffeomorphic to the total space of $U(1)$-bundle associated
with $L^{-k}$. If $S$ admits a regular Sasaki-Einstein metric, i.e. if
the underlying toric Fano manifold $M$ admits a K\"ahler-Einstein metric then the above soliton extends smoothly to the zero section.
\end{theorem}

In the set-up of the previous section \ref{linebundle}
 we next consider shrinking solitons on $L^{-k}$ with $0 < k < p$. Recall that $\lambda =  -1$ and that $a = \kappa/2 - 1 = p/k -1$ in this case.
 By the same argument as in the proof of Lemma 6.1 of \cite{FIK} we can prove that
an shrinking soliton on $L^{-k}$ must have $\mu > 0$ and $\nu = 0$.

From (\ref{GKRS3}), $\vph (a) = 0$ can be re-written as
\begin{equation}\label{GKRS6}
f(a,\mu) := \frac{(m+1)!}{a^m\mu^{m+2}}\left( \frac{2a^{m+1}\mu^{m+1}}{(m+1)!} + \left(2 - \frac{\kappa\mu}{m+1}\right) \sum_{j=0}^m \frac{a^j\mu^j}{j!}\right) = 0.
\end{equation}
The following lemma can be proved in the same way as Lemma 6.2 in \cite{FIK}.
\begin{lemma} \label{GKRS7}
For each $ 0 < \s < \frac\kappa2$, there exists a unique positive root $\mu$ of $f(\s, \mu) = 0$. This root satisfies $\mu > \frac{2(m+1)}\kappa$.
\end{lemma}
\begin{proof} We may write $f$ in the alternate form
\begin{equation}\label{GKRS8}
f(\s,\mu) = \frac{(m+1)!}{\s^m\mu^{m+2}}\left(\sum_{j=0}^{m+1} \frac{(2\s - \frac{\kappa j}{m+1})\s^{j-1}}{j!}\mu^j \right).
\end{equation}
Since $0 < \s < \kappa$, the coefficients
$$ 1, \cdots ,  \frac{(2\s - \frac{\kappa j}{m+1})\s^{j-1}}{j!}, \cdots ,  \frac{(2\s - \frac{\kappa (m+1)}{m+1})\s^{(m+1)-1}}{(m+1)!} = \frac{(2\s - \kappa)\s^m}{(m+1)!} $$
change sign only once, so there is at most one positive root $\mu$ of $f(\s, \mu) = 0$.
One sees from (\ref{GKRS6}) that
$f(\s, \frac {2(m + 1)}\kappa) = \s > 0$, while $f(\s, \mu) \sim (\s - \kappa)/\mu < 0$ as $\mu \to \infty$.
\end{proof}
Since there is no exponential term one can see that $\s$ exists for $-\infty < s < \infty$. It is now obvious that $ \s > 0$, $\s'  > 0$. This defines a shrinking soliton on $L^{-k}$.

The vector field $-\frac{\mu}{2t} r\frac \p{\p r}$ generates a one parameter group $\{\gamma_t\}$ of transformations such that
$$ \gamma_t^{\ast} r = (-t)^{-\frac\mu2} r.$$
Then
$$ \gamma_t^{\ast} s = \log \gamma_t^{\ast} r = -\frac\mu2 \log (-t) + s.$$
The flow $\{t\gamma_t^{\ast}\omega\}$ satisfies
\begin{eqnarray}
-t\gamma_t^{\ast} \omega &=& -t ((-t)^{-\frac\mu2}r)^{\frac 2\mu} D(((-t)^{-\frac\mu2} r)^{-\frac 2\mu}) i\p\bp s  \nonumber\\
&& - t \left(\frac 2\mu ((-t)^{-\frac\mu2}r)^{\frac 2\mu} D((-tr^{-\frac2\mu}) - \frac 4{\mu^2} D'(tr^{-\frac2\mu})\right) i\p s \wedge \bp s \nonumber \\
&\to & D(0)\left( r^{\frac 2\mu} i \p\bp s + \frac 2\mu r^{\frac 2\mu} i\p s \wedge \bp s \right) \label{limit3}
\end{eqnarray}
as $t \to 0$.
 If we put $q = \frac 1{\mu}$ then (\ref{limit3}) is equal to
 \begin{equation}
 D(0) \left(r^{2q} i \p\bp \log r + 2q r^{2q} i\p \log r \wedge i\bp \log r \right)
 = D(0) i\p\bp \left(\frac {r^{2q}}{2q}\right).
\end{equation}
This is a Ricci-flat K\"ahler cone metric with aperture, or equivalently a K\"ahler cone metric of a transversely K\"ahler-Einstein Sasaki manifold
with positve basic first Chern class.

Thus we have proved :
\begin{theorem}\label{shrinking(t<0)} Let $M$ be a Fano manifold, and $L \to M$ be a positive line bundle with
$L^{-p} = K_M$, $p \in \bfZ^+$.
Suppose that the $U(1)$-bundle of $K_M$ admits a possibly irregular Sasaki-Einstein metric
whose cone $C(S)$ is biholomorphic to the total space of $K_M$ minus the zero section.
For $0 < k < p$, $L^{-k}$ minus the zero section admits a gradient shrinking soliton $g(t)$ for $-\infty < t < 0$
such that $g(t)$ converges as $t \to 0$ to a Ricci-flat K\"ahler cone metric with aperture, or equivalently a K\"ahler cone metric
over a transversely K\"ahler-Einstein Sasaki manifold.
Here the K\"ahler cone manifold is biholomorphic to $L^{-k}$ minus the zero section and the
transversely K\"ahler-Einstein Sasaki manifold is diffeomorphic to the total space of $U(1)$-bundle associated
with $L^{-k}$. If $S$ admits a regular Sasaki-Einstein metric, i.e. if
the underlying toric Fano manifold $M$ admits a K\"ahler-Einstein metric then the soliton extends smoothly to the zero section.
\end{theorem}

\begin{proof}[\it{Proof of Theorem \ref{mainthm1}}  ]
 By Theorem \ref{expandingcao} we have the expanding soliton on $(L^{-k} - Z)$ and the corresponding K\"ahler-Ricci flow
 on $(L^{-k} - Z) \times (0, \infty)$.
 By Theorem \ref{shrinking(t<0)}  we also have the shrinking soliton on $(L^{-k} - Z)$ and the corresponding K\"ahler-Ricci flow
 on $(L^{-k} - Z) \times (-\infty, 0)$.
By adjusting the solitons by homothety
so that $E(0) = D(0)$ we get a smooth soliton on $(L^{-k} - Z) \times (-\infty, \infty)$.
If $S$ admits a regular Sasaki-Einstein structure then the shrinking soliton extends smoothly to the zero section as stated in
Theorem \ref{shrinking(t<0)}. This completes the proof of Theorem \ref{mainthm1}
\end{proof}

\section{Complete solitons in the cone of compact $\eta$-Einstein Sasaki manifolds}

Let us define a {\it gradient scalar soliton} to be a K\"ahler metric $g$ such that the scalar curvature $\mathcal S$ satisfy
\begin{equation}\label{GSS}
{\mathcal S} - c + \Delta Q = 0
\end{equation}
where $c$ is a constant and $Q$ is a smooth function whose gradient vector field of $Q$ is the real part of a holomorphic
vector field.
The gradient scalar solitons are also called generalized quasi-Einstein metrics (\cite{Guan95}, \cite{Maschler-Tonn09}).
We wish to find gradient scalar solitons using Calabi's ansatz on the cone of Sasaki manifold with transverse K\"ahler-Einstein
structure (or equivalently $\eta$-Eintstein Sasaki manifold), and with this purpose we go back to the beginning of section \ref{SRFKC}.
Let $\omega_{\vph}$ be the K\"ahler metric on $C(S)$ defined by Calabi's ansatz, expressed as (\ref{omega1}). In this section we require
$(a,b) = (1,\infty)$. Thus we require
$$ \s - 1 = F^{\prime}(S) > 0$$
and
$$ \vph(\s) = F^{\prime\prime} > 0.$$
Let $u(\s)$ be a smooth function of $\s$. Then
\begin{eqnarray}\label{C20}
dd^c\, u(\s) &=& d\,(u'(\s)\, \frac{d\s}{ds}\, d^c s)\\
&=& u'(\s) \varphi(\s)dd^c s + (u'\varphi)' \varphi ds \wedge d^cs \nonumber\\
&=& u'(\s) \varphi(\s)dd^c s + \frac 1{\varphi}(u'\varphi)'  d\s \wedge d^c\s. \nonumber
\end{eqnarray}
Taking wedge product of this with
\begin{equation}\label{C21}
\omega_{\varphi}^m = \s^m(\omega^T)^m + m\s^{m-1}\wedge \varphi^{-1}d\tau
\wedge d^c\s
\end{equation}
and comparing it with
\begin{equation}\label{C22}
\omega_{\varphi}^{m+1} = \s^m(m+1) \varphi^{-1}\,d \s
\wedge d^c \s \wedge(\omega^T)^m.
\end{equation}
we obtain the Laplacian $\Delta_{\vph}$ with respect to $\omega_{\vph}$ is expressed as
\begin{equation}\label{C23}
\Delta_{\varphi} u = \frac m\s u' \varphi + (u'\varphi )'.
\end{equation}
From (\ref{omega1}) and (\ref{ricci}) the scalar curvature ${\mathcal S}_{\vph}$ of $\omega_\vph$ is given by
\begin{eqnarray}
{\mathcal S}_{\vph} &=& \frac{\kappa m}\s - \Delta_{\vph} \log (\s^m \vph(\s)) \label{GSS2}\\
&=& \frac{\kappa m}\s - \frac {m\vph}\s \frac d{d\s} \log \s^m \vph - \frac d{d\s} (\vph \frac d{d\s} \log \s^m \vph)\nonumber\\
&=& \frac{\kappa m}\s - \frac 1{\s^m} \frac {d^2}{d\s^2} (\s^m\vph).\nonumber
\end{eqnarray}
From (\ref{GSS2}) and Lemma \ref{grad1} the gradient scalar soliton equation is written as
\begin{eqnarray*}
\mathcal S_{\vph} -c &=& \frac {\kappa m}\s - c - \frac 1{\s^m} \frac {d^2}{d\s^2} (\s^m\vph)\\
&=& \Delta_\vph ( -\mu \s)\\
&=& - \mu (\frac {m\vph}\s + \vph^{\prime}) \\
&=& - \frac \mu {\s^m} \frac {d}{d\s} (\s^m\vph).
\end{eqnarray*}
Namely the gradient scalar soliton equation is
\begin{equation}\label{GSS3}
(\s^m\vph)^{\prime\prime} - \mu(\s^m\vph)^{\prime} = m\kappa \s^{m-1}  - c \s^m.
\end{equation}
Integrating this we obtain
\begin{equation}\label{GSS4}
(\s^m\vph)^{\prime} - \mu\s^m\vph =
\kappa \s^{m} - \frac c{m+1} \s^{m+1} + c_1.
\end{equation}
Applying (\ref{ode}) the solution to
$$ y^{\prime} - \mu y = \kappa x^m - \frac c{m+1} x + c_1$$
is given by
\begin{eqnarray}
y &=& e^{\mu x} ( - \kappa \sum_{j=0}^m \frac {x^{m-j}}{\mu^{j+1}} e^{-\mu x} \frac {m!}{(m-j)!} +
\frac c{\mu(m+1)} x^{m+1} e^{-\mu x} \label{GSS5}\\
& &  \qquad -  \frac c \mu \sum_{j=0}^m \frac{x^{m-j}}{\mu^{j+1}} e^{-\mu x}
-  \frac {c_1}\mu e^{-\mu x} + c_2). \nonumber\\
&=& - (\kappa + \frac c\mu) \sum_{j=0}^m \frac {m!}{(m-j)!} \frac {x^{m-j}}{\mu^{j+1}} + \frac c{\mu(m+1)} x^{m+1}
- \frac {c_1}\mu + c_2e^{\mu x}\nonumber
\end{eqnarray}
Substituting $y = \s^m\vph$ and $x = \s$ into (\ref{GSS5}) we obtain the solution $\vph(\s)$ as
\begin{equation}\label{GSS6}
\vph(\s) = - (\kappa + \frac c\mu)\sum_{j=0}^m  \frac {m!}{(m-j)!} \frac {\s^{-j}}{\mu^{j+1}} +
 \frac c{\mu(m+1)} \s
- \frac {c_1}\mu \s^{-m}+ c_2e^{\mu \s}\s^{-m}
\end{equation}

In order for the solution to be complete near $\s = 1$ we need only have
$$ \vph(1) = \vph^{\prime}(1) = 0$$
by Proposition \ref{C29}. Then it follows from (\ref{GSS4}) that
\begin{equation}\label{GSS7}
 c_1 = - \kappa + \frac{c}{m+1}.
\end{equation}
Substituting (\ref{GSS7}) into (\ref{GSS4}) we get
\begin{equation}\label{GSS8}
( \s^m\vph)^{\prime} - \mu\s^m\vph = \s^m(\kappa - \frac c{m+1}\s) - \kappa + \frac c{m+1}.
 \end{equation}
 The constant $c_2$ is determined by $\vph(1) = 0$ using (\ref{GSS6}) and (\ref{GSS7}), and is given by
 \begin{equation}\label{GSS9}
 c_2 = e^{-\mu}(\kappa \sum_{j=1}^m \frac {m!}{(m-j)!} \frac {1}{\mu^{j+1}} + c \sum_{j=0}^m \frac {m!}{(m-j)!} \frac {1}{\mu^{j+2}} ).
 \end{equation}

\begin{theorem}\label{GSS10}  Let $S$ be a compact Sasaki manifold with transversely K\"ahler-Einstein metric
with $\mathrm{Ric}^T = \kappa \omega^T$, in other words $S$ is a compact $\eta$-Einstein Sasaki manifold.
Consider Calabi's ansatz (\ref{GSS3}) for the gradient scalar soliton equation (\ref{GSS}). Suppose $\kappa - \frac c{m+1} \ge 0$,
$c < 0$ and $\mu < 0$. Then there exists a solution $\vph(\s)$ giving a complete gradient scalar soliton in the cone $C(S)$.
\end{theorem}
\begin{proof}  With the constants $c_1$ and $c_2$ given by (\ref{GSS7}) and (\ref{GSS9}) we have $\vph(1) = \vph^{\prime}(1) = 0$.
We first show that for $\s > 1$ we have $\vph > 0$. Since $c < 0$ and $\s >1$ we have from (\ref{GSS8})
\begin{eqnarray}
( \s^m\vph)^{\prime} - \mu\s^m\vph &\ge& (\kappa - \frac c{m+1}\s) - \kappa + \frac c{m+1} \nonumber \\
&=& - \frac c{m+1} (\s -1) > 0. \label{GSS11}
 \end{eqnarray}
 This shows that $\vph(\s)$ can not be nonpositive for $\s > 1$.

 Thus the K\"ahler form $\omega_{\vph}$ of Calabi's ansatz (\ref{omega1}) exists for all $\s > 1$.
 We have $\vph(1) = \vph^{\prime}(1) = 0$ and $\vph(\s)$ is linear growth when $\s \to \infty$. Hence
 this metric is complete by Proposition \ref{C29}. This completes the proof of Theorem \ref{GSS10}.
\end{proof}

\begin{proof}[Proof of Theorem \ref{mainthm2}]
Comparing (\ref{GSS4}) with (\ref{GKRS2}) we see that a gradient scalar soliton is a gradient Ricci soliton
if and only if
$$ c_1 = 0\ \  \mathrm{and}\ \ 2\lambda = - \frac c{m+1}. $$
This equivalent to
\begin{equation}\label{GSS12}
c = (m+1)\kappa = -2\lambda (m+1)
\end{equation}
The assumption of Theorem \ref{GSS10} is satisfied if $\kappa < 0$ and $\mu < 0$.
But $\kappa <0$ is assumed in Theorem \ref{mainthm2} and the choice of $\mu$ is arbitrary and
we may take $\mu < 0$. This completes the proof of Theorem \ref{mainthm2}.
\end{proof}

\end{document}